\ifx\shlhetal\undefinedcontrolsequence\let\shlhetal\relax\fi

\documentstyle[12pt]{article}
\input{mssymb}

\newtheorem{theorem}{Theorem}

\newtheorem{claim}[theorem]{Claim}
\newtheorem{example}[theorem]{Example}
\newtheorem{proof}[theorem]{PROOF}
\newtheorem{conclusion}[theorem]{Conclusion}
\newtheorem{continuation}[theorem]{Continuation}
\newtheorem{concluding}[theorem]{Concluding Remarks}
\newtheorem{observ}[theorem]{Observation}
\newtheorem{definition}[theorem]{Definition}
\def\uhr{\upharpoonright}
\def\implies{\Rightarrow}
\def\vare{\varepsilon}
\def\betha{\beta}
\def\lseq{\langle}
\def\rseq{\rangle}

\def\b{\beta}
\def\a{\alpha}
\def\l{\lambda}
\def\z{\zeta}
\def\ov{\overline}
\def\lng{\langle}
\def\rng{\rangle}
\def\rest{\mathord\restriction}

\def\uni{{\cal I}}
\def\cov{{\rm cov}}
\def\eqdf{\buildrel {{\rm def}} \over =}
\def\LR{\Leftrightarrow}
\def\triangll{\triangleleft}
\def\om{\omega}
\def\Om{\Omega}
\def\al{\aleph}
\def\imply{\Rightarrow}

\title{Cardinalities of topologies with small base\footnote{Partially
supported by  The Basic research Fund, Israeli Academy of Sciences.
Publication no. 454A done 8/1991, 3-4/1993. I thank Andrzej Roslanowski for
proofreading, pointing out gaps and rewriting a part more clearly.} 
}
\author{{\bf Saharon Shelah}\\
Department of Mathematics\\
The Hebrew University, Jerusalem, Israel\\
and\\
Department of Mathematics\\
Rutgers University , New Brunswick N.J. USA
}

\date{ August 1991\\
revised March 1993\\
last revision April 20, 1993\\
May 31, 1993
}
\begin{document}
\maketitle

\begin{abstract}
Let $T$ be the family of open subsets of a topological space (not
necessarily Hausdorff or even $T_0$). We prove that if $T$ has a base of
cardinality $\leq \mu$, $\l\leq \mu<2^\l$, $\l$ strong limit of cofinality
$\aleph_0$, then T has cardinality $\leq \mu$ or $\geq 2^\l$.  This is our
main conclusion (21). In Theorem~\ref{third} we prove it under some set
theoretic assumption, which is clear when $\lambda=\mu$; then we eliminate
the assumption by a theorem on pcf from [Sh 460] motivated originally by
this. Next we prove that the simplest examples are the basic ones; they
occur in every example (for $\lambda=\aleph_0$ this fulfill a promise from
[Sh 454]). The main result for the case $\lambda=\aleph_0$ was proved in [Sh
454].
\end{abstract}
\bigskip
\centerline{ *$\qquad$ *$\qquad$ *}
\bigskip

Why does we deal with $\lambda$ strong limit of cofinality $\aleph_0$?
Essentially as other cases are closed.
\begin{example} \label{two}
If I is a linear order of cardinality $\mu$ with $\lambda$ Dedekind cuts
{\rm then} there is a topology $T$ of cardinality $\l > \mu$ with a base $B$
of cardinality $\mu$.
\end{example}

CONSTRUCTION: Let $B$ be $ \{ [-\infty , x)_I : x\in I \} $ where $[-\infty
,x )_I =\{ y \in I : I \models y < x\}$\hfill$\Box_{\ref{two}}$

Remarks: as it is well known, if $ \mu =\mu^{<\mu},\mu < \l \leq \chi =
\chi^{\mu} $ then there is a $\mu^+ $-c.c. $\mu$ -complete forcing notion Q
, of cardinality $\chi$ such that in $V^Q$ we have $2^\mu=\chi$, there is a
$\l$-tree with exactly $\mu {\ } \l$-branches (and $\leq\mu$ other branches)
hence a linear order of cardinality $\mu$ with exactly $\l$ Dedekind cuts.
As possibly $\l^{\aleph_0} >\l $, this limits possible generalizations of
our main Theorem. Also there are results guaranteeing the existence of such
trees and linear orders, e.g. if $\mu$ is strong limit singular of
uncountable cofinality, $\mu<\l\le 2^\mu$ (see [Sh 262], [Sh 355, 3.5 $+$\S
5]) and more (see [Sh 430]).

So we naturally concentrate on strong limit cardinals of countable
cofinality.  We do not try to ``save'' in the natural numbers like $n(*)+6$
used during the proof.

\begin{theorem} \label{third} {\rm (}{\bf Main}{\rm )}
Assume 
\begin{enumerate}
\item[(a)] $\lambda_n$ for 
$n<\omega$ are regular or finite cardinals, 
$2^{\lambda_n} <\lambda_{n+1}$  and 
$\lambda = \Sigma_{n<\omega} \lambda_n(\ge \aleph_0)$.

\item[(b)] $\l=\Sigma_{n<\omega}\mu_n$ (even $\mu_{n+1} \geq \l_n $) 
and $ {\beth}_3(\mu_n)<\lambda_n$ , $\lambda \leq \mu < \lambda^{\aleph_0}
(=2^\lambda)$ and cov $(\mu, \lambda^+_n , \lambda^+_n , \mu_n^+) \leq \mu$
(see Definition below, trivial when $\lambda=\aleph_0$ and easy when
$\mu=\lambda$)

\item[(c) ]Let $T$ be the family  of open subsets of a topological
space ( not necessarily Hausdorff or even $T_0$ ), and suppose that $T$ has
a base $B$ of cardinality $\le\mu$ (i.e. $B$ is a subset of $T$ which is
closed under finite intersections, and the sets in $T$ are the unions of
subfamilies of $B$).
\end{enumerate}

\underline{Then}
\begin{enumerate}
\item  The cardinality of $T$ is either at least $\lambda^{\aleph_0}
(=2^\lambda)$ or at most $\mu$.

\item In fact, if $|T|>\mu$ then for some set $X_0$ of $\lambda$ points,
$\{U\cap X_0:U\in T\}$ has cardinality $2^\lambda$. Moreover, for some
$B'\subseteq T$ of cardinality $\lambda$, $\{X_0\cap U: U \mbox{ is the
union of a subfamily of } B'\}$ has cardinality $2^\lambda$.
\end{enumerate}
\end{theorem}

\begin{definition} ([Sh 355, 5.1])
$\cov(\mu,\l^+,\l^+,\kappa)=\min\{|P|:P$ a family of subsets of $\mu$ each
of cardinality $\le \l$, such that if $a\subseteq \mu$, $|a|\le\l$ then for
some $\a<\kappa$ and $a_i\in P$ (for $i<\a$) we have $a\subseteq
\bigcup_{i<\a}a_i\}$
\end{definition}

PROOF: Suppose we have a counterexample $T$ to \ref{third}(2) (as
\ref{third}(1) follows from \ref{third}(2)) 
with a base $B$ and let $\Omega$ be the set of points of the space, so wlog
$\lambda \leq \mu = |B| < 2^\lambda$.  Our result, as explained in the
abstract, for the case $\l=\aleph_0$ was proved in [Sh 454], and see
background there; the proof as written here applies to this case too but we
usually do not mention when things trivialize for the case $\l=\al_0$; wlog
$\Omega =\bigcup B$, $\emptyset\in B$ and $B$ is closed under finite
intersections and unions. So $T$ is the set of all unions of subfamilies of
$B$.

We prove first that:

\begin{observ}\label{familyr}
For each $n$ there is a family $R$ of cardinality $\le\mu$ of partial
functions from $\l_n$ to $\mu$ such that: {\bf for every} function $f$ from
$\l_n$ to $\mu$ {\bf there is} a partition $\lng r_\zeta|\zeta<\mu_n\rng$ of
$\l_n$ (i.e. pairwise disjoint subsets of $\l_n$ with union $\l_n$) for
which $\mathop\land\limits_{\zeta<\mu_n}$ $f\uhr r_\zeta\in R$.
\end{observ}

PROOF: By assumption (b) and $2^{\l_n}<\l\le\mu$ and $\lambda$ is strong
limit of cofinality $\aleph_0\leq\mu_n$.\hfill$\Box_{\ref{familyr}}$

\begin{claim} \label{lem:an3}
Assume $Z^*$ is a subset of $\Omega$ of cardinality at most $\mu$ and $T'$
is a subfamily of T satisfying
\[\mbox{\bf (*) }\ \ \ \ \ \ (\forall U_1,U_2\in T')[U_1=U_2\ \iff\
U_1\cap Z^* = U_2\cap Z^*],\] 
$|T'| > \mu$ and $n<\omega$ .

Then we can find a subset $Z$ of $Z^*$ of cardinality $\mu_n$, subsets
$Z_{\a} $ of $Z$ and members $U_\alpha$ of $T$ and subfamilies $T_\alpha$ of
$T'$ of cardinality $>\mu$ for $\alpha < \mu_n$ such that:
\begin{enumerate}
\item[(a)] the  sets $Z_\alpha$  for  $\alpha < \mu_n$ 
are pairwise distinct

\item[(b)] for $\alpha < \mu_n$  and  $V \in T'$ we have: $V \in T_{\alpha}$
\ iff \ $V \cap Z = Z_{\alpha}\subseteq U_\a\subseteq V$.
\end{enumerate}
\end{claim}

PROOF: We shall use {\bf (*)} freely. Define an equivalence relation $E$ on
$Z^*$: 
\[x \mbox{E} y \mbox{\ \ iff\ \ } |\{U\in T^{\prime}:x\in U\Leftrightarrow
y\not\in U \}|\leq\mu\] 
(check that E is indeed an equivalence relation).

Let $Z^{\otimes} \subseteq Z^*$ be a set of representatives. Now for
$V \in T^{\prime} $ we have: 

\centerline{$(*)\ \ \ \ \ \ \ \ \ \{ U \in T^{\prime}:U \cap
Z^{\otimes}=V \cap Z^{\otimes} \} \subseteq $}

\centerline{$\subseteq\bigcup _{xEz,\{x,z\}\subseteq Z^*} \{ U \in 
T^{\prime}:z \in U\equiv x\not\in U$ but $U \cap Z^{\otimes}= V \cap
Z^{\otimes}\}\cup\{V^*\}$}

where
$V^*=\{y\in Z^*$: for the $x\in Z^\otimes$ such that $yEx$ we have
$x\in V\}$. 
\medskip

\noindent[Why? assume $U$ is in the left side i.e. $U\in T'$ and $U\cap
Z^\otimes= V\cap Z^\otimes$; now we shall prove that $U$ is in the right
side; if $U=V^*$ this is straight, otherwise for some $x\in Z^*$, $x\in
U\equiv x\not\in V^*$; as $Z^\otimes$ is a set of representities for $E$ for
some $z\in Z^\otimes$, we have $zEx$ so by the definition of $V^*$, $x\in
V^*\iff z\in V$. But as $U\cap Z^{\otimes}=V\cap Z^{\otimes}$ we have $z\in
V\iff z\in U$. Together $x\in U\iff z\notin U$ and we are done.]
\medskip

Now the right side of $(*)$ is the union of $\leq|Z^*|^2$ sets, each of
cardinality $\leq\mu$ (by the definition of $xEz$). Hence the left side in
$(*)$ has cardinality $\leq |Z^*|^2\times\mu\leq\mu$. Let
$\{V_i:i<i^*\}\subseteq T'$ be maximal such that: $V_i\cap Z^\otimes$ are
pairwise distinct and $V_i\in T'$. So clearly
$|T'|=|\bigcup\limits_{i<i^*}\{U\in T': U\cap Z^\otimes= V_i\cap Z^\otimes
\}|\le\sum\limits_{i<i^*}\mu=\mu|i^*|$, but $|T'|>\mu$ hence $|i^*|=
|\{U\cap Z^{\otimes} : U \in T^{\prime} \} | > \mu$.  Hence (as $\lambda$ is
strong limit) necessarily $|Z^{ \otimes}| \geq \lambda$, so we can let
$z_{\betha} \in Z^{\otimes} $ for $\betha < \lambda_{n}$ be distinct. For
$\alpha<\betha<\lambda_{n}$ we know that $\neg z_\a E z_\b$ hence for some
truth value ${\bf t } _{\alpha , \betha}$ we have $ |\{ U \in
T^{\prime}:z_{\alpha} \in U \equiv z_{\beta} \not\in U \equiv {\bf
t}_{\alpha ,\beta}\} |> \mu$.  But $B$ is a base of $T$ of cardinality $\le
\mu$, hence for some $V_{\alpha ,\beta} \in B $ the set 
$$ S_{\alpha ,\beta}=\{ U \in T^{\prime} : z_{\alpha} \in U
\equiv z_{\beta} \not\in U \equiv {\bf t}_{\alpha ,\betha ,}\mbox{ and
}\{ z_{\alpha},z_{\betha} \} \cap U \subseteq V_{\alpha,\betha}\subseteq U
\}$$ 
has cardinality $> \mu$.

Choose $U^1_{\a ,\b} \in S_{\a ,\b}$ such that $\mu<|S^1_{a,\b}|$
where 
$$S^1_{\a, \b}\eqdf\{ U \in S_{\a ,\b}:U \cap \{ z_\zeta : \zeta <
\l_n \}=U^1_{\a, \b} \cap \{ z_{\zeta}:\zeta < \l_n \} \},$$ 
note that $U^1_{\a,\b}$ exists as $2^{\l_n}<\l\le\mu<|S_{\a,\b}|$.

By observation~\ref{familyr} we can find a family $R$ of cardinality $\leq
\mu$, members of $R$ has the form $\overline{u}= \lng u_{\a} :\a \in r
\rng$ , where $r \subseteq \l_{n}$, $u_{\a }\in B$ such that for every
sequence $ \ov{u} =\lng u_{\a }: \a <\l_{n} \rng $ of members of $B$, there
is a partition $\lng r_\zeta : \zeta<\mu_n\rng$ of $\l_n$ (so $r_{\z }=r_{\z
}[\ov{u} ]\subseteq \l_{n}$ for $\zeta < \mu_n$) such that $\ov{u} \rest
r_{\z } \in R$ (remember $\emptyset\in B$).  Wlog if $\overline u^\ell=\lng
u_\a^\ell : \a\in r^\ell\rng\in R$ for $\ell=1,2$ then $\overline u=\lng
u_\a : \a \in r\rng \in R$ where $r=r^1\cup r^2$ and
$u_\a=\left\{\begin{array}{cl} u^1_\a & \a\in r^1
\\ u^2_\a & \a\in r^2\setminus r^1 \\
\end{array}
\right. $.

For each $V \in T'$ we can find $ \ov{u}[V]=\lng u_{\gamma} [V] :
\gamma <\l_{n} \rng$, such that (remember $\emptyset\in B$):
$$u_{\gamma}[V] \in B,$$ 
$$z_{\gamma} \in V \Rightarrow z_{\gamma} \in u_{\gamma} [V]\subseteq V,$$   
$$ z_{\gamma} \not\in V\Rightarrow u_{\gamma} [V]=\emptyset.$$ 
Clearly there is $U^2_{\a,\b}\in S^1_{\a,\b}$ such that:
\begin{description}
\item[$(**)$] for any finite subset $w \subseteq \mu_n$ and $ \a <
\b < \l_n$, the following family has cardinality $> \mu$:
\end{description}
$$S^2_{\a,\b,w}\eqdf\{ U \in S_{\a,\b}^1:(\forall\z<\mu_n)(r_{\z}
[\ov{u}[U]]=r_{\z} [\ov{u}[U^2_{\a,\b}]])\mbox{ and }$$
$$(\forall\zeta\in w)(\ov{u} [U] \rest r_{\z}=\ov{u}[U^2_{\a ,\b}]\rest
r_{\z})\}.$$

By the Erd\"os Rado theorem for some set $M \in [ \l_{n} ]^{\mu_n^+}$:
\begin{description}
\item[(a)] for every $\a<\b$ from $M$, ${\bf t}_{\a,\b}$ are the same

\item[(b)] for every $\a < \b \in M ,\gamma,\varepsilon \in M$ the
truth values of ``$z_{\gamma} \in V_{\a ,\b}$'', ``$z_{\gamma}
\in U^2_{\a ,\b}$'',  ``$ z_{\varepsilon} \in  
u_{\gamma} [ U^2_{\a ,\b}]$'' and the value of ``Min$\{ \z <\mu_n :
\gamma \in r_{\z} [ \overline{u} [U^2_{\a ,\b } ] ] \}$'' depend just
on the order and equalities between $\a,\b,\gamma$ and $\varepsilon$.
\end{description}

Let $M=\{\a(i) : i<\mu_n^+\}$ where $[i<j\imply\a(i)<\a(j)]$, let
${\rm t}$ be $0$ if $ i<j\Rightarrow {\bf t}_{\a(i),\a(j)}=$truth and
$1$ if $i<j\implies {\bf t}_{\a(i),\a(j)}=$false.
\bigskip

\underbar{Case 1} If $i<j<\mu_n^+$ and $\vare<i \lor \vare>j$ then
$z_{\a(\vare)}\not\in U^2_{\a(i),\a(j)}$. 

So for some $\zeta_1<\mu_n$
\begin{description}
\item[$\otimes$] for every $i<\mu_n^+$, $\z_1=\min \{\z : \a(i+t)\in
r_\z[\overline{u}[U^2_{\a(i),\a(i+1)}]]\}$.
\end{description}
We let $Z=\{z_{\a(i)} : i<\mu_n\}$, $Z_i=\{z_{\a(2i+t)}\}$,
$U_i=u_{\a(2i+t)}[U^2_{\a(2i),\a(2i+1)}]$.  Clearly $U_i\cap
Z\subseteq U^2_{\a(2i),(2i+1)}$ and $U_i\cap Z=Z_i$, lastly let
$T_i=S^2_{\a(2i),\a(2i+1),\{\z_1\}}$; now $Z, Z_i$, $U_i, T_i$ are as
required.
\bigskip

\underbar{Case 2:} If $i<j<\mu_n^+$ then
$$\vare<i\implies z_{\a(\vare)}\in U^2_{\a(i),\a(j)}$$
$$\vare>j\implies z_{\a(\vare)}\not\in U^2_{\a(i),\a(j)}$$
So for some $\z_1<\mu_n$, $\z_2<\mu_n$
\begin{description}
\item[$\otimes$ (a)] for $\vare<i<j<\mu_n^+$ 
$$\z_1=\min\{\z:\a(\vare)\in r_\z[\overline{u}[U^2_{\a(i),\a(j)}]]\}$$
\item[\ \ (b)] for $i<\mu_n^+$
$$\z_2=\min\{\z:\a(i+t)\in r_\z[\overline u[U^2_{\a(i),\a(i+1)}]]\}$$
\end{description}

Let $Z=\{z_{\a(i)} :i<\mu_n\}$, $Z_i=\{z_{\a(\vare)} : \vare<2i\}
\cup\{z_{\a(2i+t)}\}$,
$U_i=\cup\{u_{\a(\vare)}[U^2_{\a(2i),\a(2i+1)}]: \vare\le 2i+1\}$ and $T_i=
S^2_{\a(2i),\a(2i+1),\{\z_1,\z_2\}}$.
\bigskip

\underbar{Case 3:}  If $i<j<\mu_n^+$ then
$$\vare<i\implies z_{\a(\vare)}\not\in U^2_{\a(i),\a(j)}$$
$$\mu_n^+>\vare>j \implies z_{\a(\vare)}\in U^2_{\a(i),\a(j)}.$$
So for some $\z_1<\mu_n$, $\z_2<\mu_n$
\begin{description}
\item[$\otimes$ (a)] for $i<\mu_n^+$, $\z_1=\min\{\z : \a(i+t)\in
r_\z[\overline{u}[U^2_{\a(i),\a(i+1)}]]\}$
\item[\ \ (b)] for $i<j<\vare<\mu_n^+$, $\z_2=\min\{\z : \a(\vare)\in
r_\z[u[U^2_{\a(i),\a(j)}]]\}$ 
\end{description}
Let $Z=\{z_{\a(i)} : i<\mu_n\}$, $Z_i=\{z_{\a(\vare)} : \vare=2i+t$ or
$2i+1<\vare<\mu_n\}$ $U_i=\cup\{u_{\a(\vare)}[U^2_{\a(2i),\a(2i+1)}] :
\vare=2i+t$ or $2i+1<\vare<\mu_n\}$ and $T_i=S^2_{\a(2i),\a(2i+1),
\{\z_1,\z_2\}}$.
\bigskip

\underbar{Case 4:} If $i<j<\mu_n^+$ then
$$\varepsilon<i\implies z_{\a(\varepsilon)}\in U^2_{\a(i),\a(j)}$$
$$\varepsilon>j\implies z_{\a(\varepsilon)}\in U^2_{\a(i),\a(j)}$$
So for some $\z_1,\z_2, \z_3<\mu_n$
\begin{description}
\item[(a)] for $\varepsilon<i<j<\mu_n^+$, $\z_1=\min\{\z:\a(\varepsilon)\in
r_\z[\overline{u}[U^2_{\a(i),\a(j)}]]\}$
\item[(b)] for $i<\mu_n^+$, $\z_2=\min\{\z:\a(i+t)\in
r_\z[\overline{u}[U^2_{\a(i),\a(i+1)}]]\}$
\item[(c)] for $i<j<\varepsilon<\mu_n^+$, $\z_3=\min\{\z:\a(\varepsilon)\in
r_\z[\overline{u}[U^2_{\a(i),\a(j)}]]\}$
\end{description}

Let $Z=\{z_{\a(i)} : i<\mu_n\}$. $Z_i=\{z_{\a(\varepsilon)} :
\varepsilon<\mu_n$ and $\varepsilon\not=2i+1-t\}$,
$U_i=\cup\{u_{\a(\varepsilon)}[U^2_{\a(2i),\a(2i+1)}]$:
$\;\varepsilon<\mu_n, \varepsilon\not=2i+1-t\}$ and 
$T_i=S^2_{\a(2i),\a(2i+1),\{\z_1,\z_2,\z_3\}}$. 
Now in all cases we have chosen $Z,T_{\a}, U_{\a}, Z_{\a}(\a < \mu_n)$
as required thus finishing the proof of the
claim.\hfill$\Box_{\ref{lem:an3}}$
\vspace{8pt}

\begin{claim}\label{ione} If $Z^*\subseteq\Omega$, $|Z^*|\le\mu$, then 
$\{U\cap Z^* : U\in T\}$ has cardinality $\le\mu$.
\end{claim}

PROOF: Assume not. We can find $T'\subseteq T$ such that:
\begin{description}
\item[($\alpha$)\ ] for $U_1,U_2\in T'$ we have $U_1=U_2 \iff U_1\cap Z^*=
U_2\cap Z^*$. 
\item[($\beta$)\ ] $|T'|>\mu$.
\end{description}
By induction on $n$ we define $\lng T_\eta, Z^1_\eta, Z^2_\eta, U_\eta :
\eta\in\prod_{\ell<n}\mu_\ell\rng$ such that:
\begin{description}
\item[(a)] $T_\eta$ is a subset of $T'$ of cardinality $>\mu$
\item[(b)] if $\nu\triangll\eta$ then $T_\eta\subseteq T_\nu$
\item[(c)] if $\eta=\lng \rng $ then 
$T_\eta=T'$, $Z_\eta^1=Z_\eta^2=\emptyset$, $U_\eta=\emptyset$
\item[(d)] $Z^1_\eta\subseteq Z^2_\eta\subseteq Z^*$ and
$|Z^2_\eta|\leq \mu_{\lg\eta},Z_\eta^2$ 
disjoint to $\cup\{U_{\eta\uhr\ell }:\ell<\lg \eta\}$ 
\item[(e)] $U_\eta\in T$
\item[(f)] if $V\in T_\eta$ then $U_\eta\subseteq V$ and $V\cap
Z^2_\eta=U_\eta\cap Z^2_\eta=Z_\eta^1$
\item[(g)] if $\lg (\eta)=\lg(\nu)=n+1$ and $\eta\uhr n=\nu\uhr n$
then $Z^2_\eta=Z^2_\nu$ 
but
\item[(h)] if $\lg(\eta)=\lg(\nu)=n+1$, $\eta\uhr n=\nu\uhr n$ but
$\eta\not=\nu$ then $Z^1_\eta\not=Z^1_\nu$.
\end{description}

Why this is sufficient? Let $Z\stackrel{\rm df}{=}\bigcup\{Z^2_\eta:
\eta\in\bigcup_n\prod_{l<n}\mu_l\}$. It is a subset of $Z^*$ of cardinality
$\leq\lambda$. The set $B'\stackrel{\rm df}{=}
\{U_\eta:\eta\in\bigcup_{n<\omega}\prod_{l<n}\mu_l\}$ is included in $T$ and
has cardinality $\leq\lambda$. For $\eta\in\prod_n \mu_n$ we let
$U_\eta=\bigcup_{n<\om}U_{\eta\uhr n}$. Now as $U_{\eta\uhr n}\in T$ (by
clause (e)), clearly $U_\eta\in T$. Now suppose $\eta\not=\nu$ are in
$\prod_{n<\om}\mu_n$ and we shall prove that $U_\eta\cap Z\not=U_\nu\cap Z$,
as $|\prod_n\mu_n|=2^\lambda$ this suffices (giving (1) + (2) from Theorem
2).  Let $n$ be minimal such that $\eta(n)\not=\nu(n)$, so $\eta\uhr
n=\nu\uhr n$. By clause (g), $Z^2_{\eta\uhr(n+1)}=Z^2_{\nu\uhr(n+1)}$. So (by
clause (h)) $Z^1_{\eta\uhr(n+1)}, Z^1_{\nu\uhr(n+1)}$ are distinct subsets
of $Z^2_{\eta\uhr(n+1)}=Z^2_{\nu\uhr(n+1)}\subseteq Z$.  So it suffices to
show $U_\eta\cap Z^2_{\eta\uhr(n+1)}=Z^1_{\eta\uhr(n+1)}$ and $U_\nu\cap
Z^2_{\nu\uhr(n+1)}= Z^1_{\nu\uhr(n+1)}$ and by symmetry it suffices to prove
the first. Now $Z^1_{\eta\uhr(n+1)}\subseteq U_{\eta\uhr(n+1)}$ by clause
(f), hence $Z^1_{\eta\uhr(n+1)}
\subseteq U_\eta$ so it suffices to prove that $U_\eta\cap
Z^2_{\eta\uhr(n+1)}\subseteq Z^1_{\eta\uhr(n+1)}$; for this it
suffices to prove that for $\ell<\om$ 
$$(*)\ \ \ \ U_{\eta\uhr \ell}\cap Z^2_{\eta\uhr(n+1)}\subseteq
Z^1_{\eta\uhr(n+1)}.$$ 
\medskip

\underbar{Case 1:} $\ell=n+1$. This holds by clause (f).
\medskip

\underbar{Case 2:} $\ell>n+1$. Then choose any $V\in T_{\eta\uhr
\ell}$, so we know $U_{\eta\uhr \ell}\subseteq V$ (by clause (f)) and $V\in
T_{\eta\uhr(n+1)}$ (by clause (b)), and $V\cap Z^2_{\eta\uhr(n+1)}=
Z^1_{\eta\uhr(n+1)}$ (by clause (f)), together finishing.
\medskip

\underbar{Case 3:}  $\ell\le n$.
By clause (d), $Z^2_{\eta\uhr (n+1)}$ is disjoint from $U_{\eta\uhr
\ell}$.
\medskip

So we have finished to prove sufficiency, but we still have to carry the
induction. For $n=0$ try to apply (c), the main point being
$|T_{\lng\rng}|>\mu$ which holds by the choice of $T'$ (which was possible
by the assumption that the claim fails).  Suppose we have defined for $n$
and let $\eta\in\prod_{\ell<n}\mu_\ell$.  We apply claim \ref{lem:an3} with
$T_\eta$, $Z^*\setminus\bigcup\limits_{\ell< n} U_{\eta\uhr \ell}$ and $n$
here standing for $T'$, $Z^*$, $n$ there.

We get there $Z,Z_\a,T_\a,U_\a$ ($\a<\mu_n$) satisfying (a)$+$(b)
there. We choose $T_{\eta^\land\lng\a\rng}$ to be $T_\a$,
$U_{\eta^\land\lng\a\rng}$ to be $U_\a$, $Z^2_{\eta^\land\lng\a\rng}$
to be $Z$ and $Z^1_{\eta^\land\lng\a\rng}$ to be $Z_\a$. You can check
the induction hypotheses, so we have finished.\hfill$\Box_{\ref{ione}}$

\begin{definition} 
\label{small}
$X\subseteq \Omega$ is \underbar{small} if $\{X\cap U:U\in T\}$ has
cardinality $\le\mu$. The family of small $X\subseteq\Omega$ will be denoted
by ${\cal I}={\cal I}_T$ (or more exactly, ${\cal I}_{T,\Omega}$)
\end{definition}

\begin{claim} \label{lem:an5}
The family of small sets, ${\cal I}$, is a $\mu^+$-complete ideal (on
$\Om$, including all singletons of course).
\end{claim}

PROOF: Clearly $\uni $ is a family of subsets of $\Omega$, and it is
trivial to check that $X\in \uni$ and $Y\subseteq X\implies Y\in
\uni$.  So assume $X_\alpha\in \uni$ for $\alpha<\alpha(*)$,
$\alpha(*)\le\mu$ and we shall prove that $X=\bigcup_\alpha
X_\alpha\in \uni$.  Each $X_\alpha$ has a subset $Y_\alpha$ such that
\begin{description}
\item[(a)] $|Y_\alpha|\leq\mu$ and 
\item[(b)] if $V, W$ are elements of $T$ with $V\cap X_\alpha\not = W\cap
X_\alpha$ then there is some element $y \in Y_\alpha$ which is in exactly
one of $V, W$ (possible as $X_\a\in\uni$). 
\end{description}
Now if $V, W$ are elements of $T$ which differ on
$X=\bigcup_{\alpha<\alpha(*)} X_\alpha$, then they already differ on some
$X_\alpha$ and hence they differ on some $Y_\alpha$ hence on $Y \eqdf
\bigcup_{\a<\a(*)} Y_\a$.  So $|\{U\cap X: U\in T\}|=|\{U\cap Y:U\in T\}|$,
so it suffice to prove that $Y$ is small.  But $Y$ has cardinality
$\leq|\bigcup_\a Y_\a|\le\sum_\a|Y_\a|\le\mu\times\mu=\mu$; so claim
\ref{ione} implies that $Y$ is small and hence $X$ is small.\hfill
$\Box_{\ref{lem:an5}}$
\vspace{10pt}

\begin{conclusion} \label{hyp2} Wlog  
$card(\Omega) = \mu^+ $
\end{conclusion}

PROOF: As obviously $\{x\}\in\uni$ for $x\in\Omega$, by claim
\ref{lem:an5} we know $|\Om|>\mu$.  Let $T'\subseteq T$ be of
cardinality $\mu^+$ and let $\Omega'\subseteq \Omega$ be of cardinality
$\mu^+$ such that: if $U\not=V$ are from $T'$ then $U\cap
\Omega'\not= V\cap\Omega'$. Let $T''$ be $\{U\cap\Omega' : U\in T\}$
and $B'=\{U\cap\Omega': U\in B\}$.  Now $T'',B',\Omega'$ are also a
counterexample to the main theorem and satisfies the additional
demand.
\hfill$\Box_{\ref{hyp2}}$

\begin{claim}\label{nstar} 
Wlog for some $n(*)$, for no $Z\subseteq \Omega$ of cardinality $\mu_{n(*)}$
and $U_\a$, $T_\a$, $Z_\a (\a<\mu_{n(*)})$ does the conclusion of claim
\ref{lem:an3} (with $\Om$, $T$ here standing for $Z^*$, $T'$ there) holds.
\end{claim}

PROOF: Repeat the proof of claim \ref{ione}.  I.e. we let $Z^*
\eqdf\Om$, and add the demand 
$$T_\eta=\{U\in T: U_{\eta\rest l}\subseteq U\mbox{ and } U\cap
Z^2_{\eta\rest l}\subseteq U_{\eta\rest l}\mbox{ for }
l<\lg\eta\}.\leqno(i)$$ 
The only change is in the end of the paragraph before
the last one where we have used claim \ref{lem:an3}, now instead we
say that if we fail then for our $n$, replacing $T,\Om$ by $T_\eta$,
$Z^*\setminus\bigcup_{\ell < n}U_{\eta \uhr \ell}$ resp.  gives the
desired conclusion (note $T_\eta$ has a basis of cardinality $\leq\mu$: 
\[B_\eta\eqdf\{U\cup\bigcup_{l<\lg\eta}U_{\eta\rest l}: U\in B\mbox{ and }
U\cap Z^2_{\eta\rest l}\subseteq U_{\eta\rest l}\mbox{ for } l<\lg\eta\}\]
which is included in $T_\eta$).\hfill$\Box_{\ref{nstar}}$

\begin{observ} \label{lem:an4}
Suppose $\l $ is strong limit of cofinality $\aleph_0$, $I$ is a linear
order of cardinality $\le \mu$, $\l \leq \mu < \l^{\aleph_0} , $ and $I$ has $
> \mu $ Dedekind cuts, then it has $\geq \mu^{\aleph_0} (=\l^{\aleph_0}) $
Dedekind cuts.
\end{observ}

\noindent\underline{Remark}: This observation does not relay on the
assumptions of Theorem 2.
\bigskip

PROOF: We define by induction on $\a$ when does $rk_I (x,y) = \a $ for $x
<y$ in $I$.

\underline{for $\a = 0$}{\ } $rk_I (x ,y)=\a$ iff $(x ,y)_I=\{z\in
I:x<z<y\}$ has cardinality $< \l$

\underline{for $\a > 0$} {\ }$rk_I (x ,y) =\a $ if: for $\b <\a , \neg
[rk_I(x,y)=\b]$ but \underline{for any} $(x_i,y_i)$ $(i<\l)$, pairwise
disjoint subintervals of $(x,y)$,
\underbar{there is} $i$ such that $ \bigvee_{\b <\a} rk_I (x_i ,y_i )=\b $

\begin{description}
\item[$(*)_1$] Note that by thinning the family, without loss of
generality, $[x_i ,y_i]$ are pairwise disjoint,
\end{description}
[why? e.g. as for every $j$ the set $\{i: [x_i,y_i]\cap
[x_j,y_j]\neq\emptyset\}$ has at most three members].
\begin{description}
\item[$(*)_2$] for $\a>0$ and $x<y$ from $I$, $rk_I(x,y)=\a$ iff for
$\b<\a$, $\neg[rk_I(x,y)=\b]$ and for some $\l^{\prime} < \l$ for any
$(x_i ,y_i)$ $(i <\l^{\prime} )$, pairwise disjoint subintervals of
$(x,y)$ there are $i<\l^{\prime}$ and $\b < \a$ such that $rk_I
(x_i,y_i) = \b$
\end{description}
\noindent [Why? the demand in $(*)_2$ certainly implies the demand in the
definition, for the other direction assume that the definition holds but the
demand in $(*)_2$ fails, and we shall derive a contradiction.  So for each
$n < \om$ there are pairwise disjoint subintervals $(x^n_i , y^n_i )$ of
$(x,y)$, for $i < \l_n$ such that $\neg [rk_I ( x^n_i , y^n_i ) = \b ]$
(when $\b < \a$ and $i < \l_n$).  As we can successively replace $\{(x^n_ i
, y^n_i ) : i <\l_n \}$ by any subfamily of the same cardinality (when the
$\l_n$'s are finite - by a subfamily of cardinality $\l_{n-1}$) wlog: for
each $n$, all members of $\{ x^n_i : i < \l_n \}$ realize the same Dedekind
cut of $\{ x^m_j , y^m_j : m<n,j < \l_m \}$ and similarly for all members of
$\{ y^n_i : i < \l_n \}$.  So for $m < n , i< \l_n$, the interval $(x^n_ i ,
y^n_i )$ cannot contain a point from $\{x^m_j , y^m_j : j < \l_m \}$ (as
then the same occurs for all such $i$'s, for the same point contradicting
the ``pairwise disjoint'') so either our interval $(x^n_ i , y^n_i)$ is
disjoint to all the intervals $(x^m_j , y^m_j)$ for $j< \l_m$ or it is
contained in one of the intervals $(x^m_j , y^m_j)$; as $j$ does not depend
on $i$ we denote it by $j(m,n)$; if $\lambda=\aleph_0$, by the Ramsey
theorem wlog for $m<n$ $j(m,n)$ does not depend on $n$; now the family
$\{(x^m_i , y^m_i ): m <
\om$, $i <\l_m$ and for every $n < \om$ which is $> m$ we have $i \not=
j(m,n) \}$ contradicts the definition]
\bigskip

\noindent If $rk_I (x ,y) $ is not equal to any ordinal let it be $\infty$.
Let $\a^*=\sup\{rk_I(x,y)+1 : x<y$ in $I$ and $rk_I(x,y)<\infty\}$.  Clearly
$rk_I(x ,y)\in\a^*\cup\{\infty\}$ for every $x <y $ in $I$ (and in fact
$\a^*<\mu^+ $ ). As we can add to $I$ the first and the last elements it
suffices to prove:

(A) if $rk_I (x ,y) = \a < \infty $ then $ (x ,y)_I $
has $\leq \mu $ Dedekind cuts and 

(B) if $rk_I (x ,y) =\infty $ then it has $\geq \l^{\aleph_0} $
Dedekind cuts

\noindent (B) is straightforward.

\noindent{\em Proof of }(A): We prove this by induction on $\a$. If $\a$ is
zero this is trivial. So assume that $\a > 0$, hence by $(*)_2$ for some
$\l'<\l$ there are no pairwise disjoint subintervals $(x_i , y_i )$ for $i
<\l'$ such that $\b < \a$ implies $\neg[rk_I(x_i,y_i)=\b]$.  Let $J$ be the
completion of $I$, so each member of $J\setminus I$ realizes on $I$ a
Dedekind cut with no last element in the lower half and no first element in
the upper half, and $|J| > \mu \ge |I|$. Let $J^+ \eqdf \{ z \in J: z
\not\in I$ and if $x \in I$, $y \in I$ and $x <_J z <_J y$ and $\b < \a$
then $\neg [rk_I (x,y ) =\b] \}$.  By the induction hypothesis, easily $|J
\setminus J^+ | \le \mu$ hence the cardinality of $J^+$ is $> \mu$. By
Erd\"os-Rado theorem, (remembering $\l$ is strong limit and $\l' < \l$)
there is a monotonic (by $<_J$) sequence $\langle z_i:i< \l' \rangle$ of
members of $J^+$; by symmetry wlog $\langle z_i : i< \l'
\rangle$ is $<_J$ -increasing. Now for each $i < \l'$ as $z_i < _J
z_{i+1}$ both in $J^+$ neccessarily there is a member $x_i$ of $I$ such that
$z_i <_J x_i <_J z_{i+1}$.  So $x_i <_J z_{i+1} <_J x_{i+1}$ and $x_i
\in I$, $x_{i+1} \in I$ and $z_{i+1} \in J^+$ hence by the definition
of $J^+$ we know that for no $\b < \a$ is $rk_I (x_i , x_{i+1} ) =
\b$. So finally the family $\{ (x_i , x_{i+1} ) : i < \l' \}$ of
subintervals of $(x,y)$ gives the desired contradiction to $(*)_2$.
\hfill$\Box_{\ref{lem:an4}}$  

\begin{definition}\label{12B}
We define an equivalence relation $E$ on $\Omega$: $xEy$ iff $\{U\in T: x\in
U\equiv y\not\in U\}$ has cardinality $\leq \mu$. 
\end{definition}

\begin{conclusion} \label{hyp3}
\begin{description}
\item[(0)] The equivalence relation $E$ has $< \l_{n(*)}<\l$
equivalence classes (for some $n(*)<\om$, which wlog is as required in
claim \ref{nstar} too).
\item[(1)] wlog  for each $x \in \Omega$ one of the following sets has
cardinality $\leq \mu$ : 

(a) $\{ U \in T : x \in U\}$

(b) $\{ U \in T : x \notin U\}$ 

\item[(2)] wlog for all $x \in \Omega$  we get the same case above, in fact
it is case (b). 

\item[(3)] wlog for any two distinct  members $x,y$
of $\Omega$  for some $U \in B$  we have
$ x \in U$  iff $ y  \notin U$.
\end{description}
\end{conclusion}

PROOF: (0) By claim \ref{nstar} and the proof of claim \ref{lem:an3}
(if $E$ has $\ge \l$ equivalence classes we can repeat the 
proof of claim \ref{lem:an3} and get contradiciton to claim \ref{nstar}).

(1), (2), (3) Let $\lng X_\z:\z<\z^*\rng$ list the $E$-equivalence
classes, so $\z^*<\l_{n(*)}$.  As $\Om\not\in\uni$, and $\uni$ is
$\mu^+$-complete (claim \ref{lem:an5}) for some $\z$,
$X_\z\not\in\uni$. Let $\Om'=X_\z$, $T'=\{U\cap\Om':U\in T\}$,
$B'=\{U\cap\Om':U\in B\}$; so $\Om',B',T'$ has all the properties we
attribute to $\Om,B,T$ and in addition now $E$ has one equivalence
class. So we assume this.

Fix any $x_0\in\Om$, let $B^0=\{U\in B:x_0\not\in U\}$, $T^0=\{U\in
T:x_0\not\in U\}\cup\{\Omega\}$, $B^1=\{U\in B:x_0\in U\}$, $T^1=\{U\in
T:x_0\in U\}\cup\{\emptyset\}$. For some $\ell\in\{0,1\}$, $|T^\ell|>\mu$,
and then $\Om,B^\ell,T^\ell$ satisfies the earlier requirements and the
demands in (1) and (2). For (3) define an equivalence relation $E'$ on
$\Om$: $xE'y$ iff $(\forall U\in B)[x\in U\equiv y\in U]$, let
$\Om'\subseteq\Om$ be a set of representatives, $B'=\{U\cap\Omega':U\in B\}$
and finish as before. The only thing that is left is the second phrase in
(2). But if it fails then for every $U\in T\setminus\{\emptyset\}$ choose a
nonempty subset $V[U]$ from $B$. As the number of possible $V[U]$ is
$\leq|B|\leq\mu$, for some $V\in B\setminus\{\emptyset\}$, for $>\mu$
members $U$ of $T$, $V=V[U]$ and hence $V\subseteq U$. Choose $x\in V$; so
for $x$ clause (a) of (2) fails and hence for all $y\in\Omega$ clause (b) of
(2) holds, as required. \hfill$\Box_{\ref{hyp3}}$


\begin{proof}
{\bf (of Theorem \ref{third} (MAIN)):}
\end{proof}
Consider for $n=n(*)$ (from claim \ref{hyp3}(0) and as in claim \ref{nstar})
the following:

$(*)$ there are an open set $V$ and a subset $Z$ of $V$ and for each
$\alpha < \lambda_n$ $Z_\a\subseteq Z$ and open subsets $V_\alpha,U_\a
$ of $V$ such that: 

(a) for $\alpha < \betha < \lambda_n$ the sets $V_\alpha \cap Z, V_\betha
\cap Z$ are distinct 

(b) $U_{\a}\cap Z=Z_{\a}$

(c) the number of sets $U \in T$ satisfying $U \cap Z = V_\alpha \cap
Z$ and $U_\a \subseteq U$ is $>\mu$

So by claim \ref{nstar} we know that this fails for $n$.

Let $\chi$ be large enough and let $\bar N=\lseq N_i : i < \mu^+
\rseq$ be an elementary chain of submodels of $(H(\chi), \in )$ of
cardinality $\mu$ (and $B,\Omega, T$ belong to $N_0$ of course) increasing
fast enough hence e.g.: if $X\in N_i$ is a small set, $U\in T$ then there is
$U'\in N_i\cap T$ with $U\cap X=U'\cap X$ (you can avoid the name
"elementary submodel " if you agree to list the closure properties actually
used; as done in [Sh 454]). For $x \in
\Omega$ let $i(x)$ be the unique $i$ such that $x$ belongs to $N_{i
+1} \backslash N_i$ or $i=-1$ if $x \in N_0$ (remember $|\Om|=\mu^+$).

\begin{definition}
 We define : $x \in \Omega$ is $\overline{N}$-pertinent if it belongs
to some small subset of $\Omega$ which belongs to $N_{i(x)}$ (and
$i(x)\geq 0$) and $\overline{N}$-impertinent otherwise.
\end{definition}

\begin{observ} \label{ob1}
$\Omega_{ip}=\{x\in\Omega: x$ is $\overline{N}$-impertinent $\}$ is
not small (see Definition~\ref{small}).
\end{observ}

PROOF: As $N_0\cap\Omega$ is small by claim \ref{lem:an5}, for some
$U^*$, $T'\eqdf\{U\in T: U\cap N_0\cap\Omega= U^*\cap N_0\cap\Omega\}$
has cardinality $>\mu$. So it suffices to prove:

$(*)\;\;U_1\not=U_2\in T'\implies
U_1\cap\Omega_{ip}\not=U_2\cap\Omega_{ip}$.

Choose $x\in (U_1\setminus U_2)\cup(U_2\setminus U_1)$ with $i(x)$
minimal.  As $U_1, U_2\in T'$, $i(x)=-1$ ( i.e. $x\in N_0$) is
impossible, so $x\in(N_{i+1}\setminus N_i)\cap\Omega$ for $i=i(x)$.
If $x\in \Omega_{ip}$ we succeed so assume not i.e. $x$ is
$\overline{N}$-pertinent, so for some small $X\in N_i$ $x\in X$. Hence
by the choice of $\bar N$: for some $U'_1,U'_2\in N_i\cap T$ we have:
$U'_1\cap X=U_1\cap X, U'_2\cap X=U_2\cap X$ so $U'_1\cap X$,
$U'_2\cap X\in N_i$ are distinct (as $x$ witness) so there is $x'\in
N_i\cap X$, $x'\in U'_1\equiv x'\not\in U'_2$; but this implies $x'\in
U_1\equiv x'\not\in U_2$, contradicting $i(x)$'s minimality.
\hfill$\Box_{\ref{ob1}}$
\bigskip

We define a binary relation $\preceq$ on $\Om_{ip}$ by: 
$$x\preceq y\LR \ {\rm for\ all }\ U\in B,\ {\rm if }\ y\in U \ {\rm
then }\ x\in U.$$

\begin{claim}
The relation $\preceq$ is clearly reflexive and transitive. It is
antisymetric [why antisymetric? by claim \ref{hyp3}(3)].
\end{claim}

\begin{observ} 
\label{lin:order}
If $J\subseteq \Om_{ip}$ is linearly ordered by 
$\preceq$ then $J$ is small.
\end{observ}

PROOF: For each $U_1,U_2\in B$ such that $U_1\cap J\not\subseteq U_2\cap J$
choose $y_{U_1,U_2}\in J\cap (U_1\setminus U_2)$. Let $I=\{y_{U_1,U_2}:
U_1,U_2\in B\ \&\ U_1\cap J\not\subseteq U_2\cap J\}$.  Clearly
$|I|\leq\mu$. We claim that $I$ is dense in $J$ (with respect to $\preceq$,
i.e. $I$ has a member in every non empty interval of $J$).  Suppose that
$x,y,z\in J$, $x\prec y\prec z$. By \ref{hyp3}(3) we find $U_1, U_2\in B$
such that $x\in U_1, y\notin U_1$, and $y\in U_2, z\notin U_2$. Consider
$y_{U_2, U_1}\in I$. Easily $x\prec y_{U_2,U_1}\prec z$. Thus if
$(x,z)\neq\emptyset$ then $(x,z)\cap I\neq\emptyset$.

Now note that each Dedekind cut of $I$ is an restriction of at most 3
Dedekind cuts of $J$ (and the restriction of a Dedekind cut of $J$ to $I$ is
a Dedekind cut of $I$). For this suppose that $Y_1, Y_2, Y_3, Y_4$ are lower
parts of distinct Dedekind cuts of $J$ with the same restriction to $I$,
wlog $Y_1\subset Y_2\subset Y_3\subset Y_4$. For $i=2,3,4$ choose $y_i\in
Y_i$ such that $Y_1\prec y_2, Y_2\prec y_3$ and $Y_3\prec y_4$. As $(y_2,
y_4)\neq\emptyset$ we find $x\in (y_2,y_4)\cap I$. Since $y_2\prec x$ we get
$x\notin Y_1$ and since $x\prec y_4$ we obtain $x\in Y_4$. Consequently $x$
distinguishes the restrictions of cuts determined by $Y_1$ and $Y_4$ to $I$.

To finish the proof of the observation apply observation \ref{lem:an4} to
$I$ (which has essentially the same number of Dedekind cuts as $J$).
\hfill$\Box_{\ref{lin:order}}$

\begin{continuation} {\bf (of the proof of theorem \ref{third})}
\end{continuation}
Now it suffices to prove that for each $x\in\Omega_{ip}, i=i(x)>0$ there is
no member $y$ of $\Omega_{ip} \cap N_i$ such that
$x,y$ are $\preceq$-incomparable.

\noindent [Why? then we can divide $\Omega_{ip}$ to $\mu$ sets such that any
two in the same part are $\preceq$-comparable contradicting 16$+$18 and 8; How?
By defining a function $h:\Omega_{ip}\longrightarrow\mu$ such that $h(x)=h(y)
\Rightarrow x \preceq y \vee y\preceq x$. We define $h\rest(\Omega_{ip}
\cap N_i)$ by induction on $i$, in the induction step let $N_{i+1}\backslash
N_i =\{ x_{i,\varepsilon} {\ } : \varepsilon <\mu \}$.  Choose $h
(x_{i,\varepsilon} )$ by induction on $\vare$: for each $\varepsilon$ there
are $\leq|\varepsilon | < \mu $ forbidden values so we can carry the
definition.]

So assume this fails, so we have: for some $x \in \Omega_{ip}, i=i(x)>0$
there is $y_0\in N_i\cap\Omega_{ip}$ which is $\preceq$-incomparable with
$x$; so there are $U_0$, $V_0\in B$ such that $x\in V_{0}$, $x\notin U_{0}$,
$y_0\in U_{0}$, $y_{0} \notin V_{0}$. Now $U^*=\bigcup\{U\in T: y_0\notin
U\}$ is in $T\cap N_i$ and $x\in U^*$ (as $V_0$ witnesses it) but by
\ref{hyp3}(2) we know that $U^*$ is small, so it contradicts
``$x\in\Omega_{ip}$''. This finishes the proof of theorem
\ref{third}.\hfill$\Box_{\ref{third}}$ 

\begin{concluding}\label{last}
Condition (b) of Theorem \ref{third} holds easily for $\mu=\l$. Still it may
look restrictive, and the author was tempted to try to eliminate it (on such
set theoretic conditions see [Sh 420,\/\S6]).  But instead of working
``honestly'' on this the author for this purpose proved (see [Sh 460]) that
it follows from ZFC, and therefore can be omitted, hence
\end{concluding}

\begin{conclusion}\label{main} {\rm (}{\bf Main}{\rm )} If $\l$ is
strong limit, cf$\l=\aleph_0$, and $T$ a topology with base $B$,
$|T|>|B|\ge\l$ then $|T|\ge2^\l$ and thew conclusion of \ref{third}(2)
holds. 
\end{conclusion}

\begin{theorem}
\label{22}
\begin{enumerate}
\item Under the assumptions of Theorem~\ref{third}, if the topology $T$ is
of the size $\geq 2^\lambda$ then there are distinct $x_\eta\in\Omega$ for
$\eta\in\bigcup_{n<\omega}\prod_{l<n}\mu_l$ such that letting
$Z=\{x_\eta:\eta\in\bigcup_{n<\omega}\prod_{l<n}\mu_l\}$ one of the
following occurs:
\begin{description}
\item[(a)\ ] there are $U_\eta\in T$ (i.e. open) for
$\eta\in\prod_{l<\omega}\mu_l$ such that:
\[U_\eta\cap Z=\{x_\nu\in Z: (\exists n<\mbox{lg\/}(\nu))(\nu\rest n=\eta\rest
n \ \&\ \nu(n)<\eta(n))\}\]
\item[(b)\ ] there are $U_\eta\in T$ for $\eta\in\prod_{l<\omega}\mu_l$ such
that:  
\[U_\eta\cap Z=\{x_\nu\in Z: (\exists n<\mbox{lg\/}(\nu))(\nu\rest n=
\eta\rest n \ \&\ \nu(n)>\eta(n))\}\]
\item[(c)\ ] there are $U_{\eta}\in T$ for
$\eta\in\prod_{l<\omega}\mu_l$ such that:
\[U_{\eta}\cap Z=\{x_\nu\in Z: \neg\nu\triangll\eta\}\]
\end{description}
\item If in addition $\lambda=\aleph_0$ then we get\\
$\oplus$\ \ \ \ \ there are distinct $x_q\in\Omega$ for $q\in{\Bbb Q}$ (the
rationals) such that for every real $r$, for some (open) set $U\in T$
\[U\cap\{x_q: q\in {\Bbb Q}\}=\{x_q: q\in {\Bbb Q}, q<r\}.\] 
\end{enumerate}
\end{theorem}

\begin{observ}
\label{obser}
Suppose that there are distinct $x_\eta\in\Omega$ (for
$\eta\in\bigcup_{n\in\omega}\prod_{l<n}\mu_l$) such that one of the
following occurs:
\begin{description}
\item[(d)\ ] there are $U_\eta\in T$ for $\eta\in\prod_{l<\omega}\mu_l$ such
that:
\[U_\eta\cap Z=\{x_\nu\in Z: \nu=\rho\hat{\ }\langle\zeta\rangle\ \&\ 
[\neg\rho\triangll\eta\mbox{ or }\rho\triangll\eta\ \&\
\eta(\mbox{lg\/}(\rho))=\zeta ]\}\] 
\item[(e)\ ] there are $U_\eta\in T$ for $\eta\in\prod_{l<\omega}\mu_l$ such
that:
\[U_\eta\cap Z=\{x_\nu\in Z: \nu=\rho\hat{\ }\langle\zeta\rangle\ \&\ 
[\neg\rho\triangll\eta\mbox{ or }\rho\triangll\eta\ 
\&\ \eta(\mbox{lg\/}(\rho))<\zeta ]\}\]
\item[(f)\ ] there are $U_\eta\in T$ for $\eta\in\prod_{l<\omega}\mu_l$ such
that:
\[U_\eta\cap Z=\{x_\nu\in Z: \nu=\rho\hat{\ }\langle\zeta\rangle\ \&\
[\neg\rho\triangll\eta\mbox{ or }\rho\triangll\eta\ 
\&\ \eta(\mbox{lg\/}(\rho))>\zeta ]\}.\]
\end{description}
Then  for some distinct $x_\nu'\in\Omega$ ($\nu\in\bigcup_{n\in\omega}$) the
clause {\bf (c)} of theorem~\ref{22} holds.
\end{observ}

PROOF Let $U_\eta$ (for $\eta\in\prod_{l\in\omega}\mu_l$) be given by one
of the clauses. For $\nu\in\prod_{l<n}\mu_l$, $n\in\omega$ let
$g(\nu)\in\prod_{l<2n}\mu_l$ be such that $g(\nu)(2l)=0$,
$g(\nu)(2l+1)=\nu(l)$ and for $\eta\in\prod_{l\in\omega}\mu_l$ let
$g(\eta)=\bigcup_{l<\omega}g(\eta\restriction l)$ (we assume that
$\mu_l<\mu_{l+1}$). Next define points $x_\nu'\in\Omega$ and open sets
$U_\eta'$ as 
\[U_\eta'=U_{g(\eta)},\ \ \ \ x_\nu'=\left\{
\begin{array}{ll}
x_{g(\nu)\hat{\ }\langle 1\rangle} & \mbox{if we are in clause {\bf (d)}}\\
x_{g(\nu)\hat{\ }\langle 0\rangle} & \mbox{if we are in clauses {\bf (e), (f)}}
\end{array}
\right.\]
Then $x_\nu'$, $U_\eta'$ examplify clause {\bf (c)} of theorem~\ref{22}
\hfill$\Box_{\ref{obser}}$


\begin{proof}
\label{subcase}
of \ref{22} for the case $\lambda=\aleph_0$
\end{proof}
It suffices to prove \ref{22}(2), as $\oplus$ implies {\bf (a)}. 
Let $\mu=\lambda^+$. By Theorem~\ref{third}(2) and \ref{main} wlog
$|\Omega|=\lambda$, $|B|\leq\lambda$. Let $\uni = \{Z\subseteq\Omega:
|\{U\cap Z: U\in T\}|<\mu\}$, again it is a proper ideal on $\Omega$ (but
not necessarily even $\aleph_1$-complete). Let $P=\{(U,V): U\subseteq V
\mbox{ are from T, } V\setminus U\notin\uni\}$.  Clearly $P\neq\emptyset$
(as $(\emptyset,\Omega)\in P$), if for every $(U_0,U_1)\in P$ there is $U$
such that $(U_0,U), (U,U_1)$ are in $P$ then we can easily get clause 
$\oplus$. So by renaming wlog  
$$(\forall V\in T)(V\in\uni\ \ \mbox{or}\ \ \Omega\setminus
V\in\uni).\leqno(*)_1$$  
We try to choose by the induction on $n<\omega$, $(x_n, U_n)$ such that
\begin{description}
\item[(a)\ ] $x_n\in U_n\in T$
\item[(b)\ ] $x_n\notin\bigcup_{l<n}U_l$
\item[(c)\ ] $U_n\in\uni$ and $x_l\notin U_n$ for $l<n$
\item[(d)\ ] $|\{V\in T: (\forall l\leq n)(x_l\notin V)\}|\geq\mu$.
\end{description}
If we succeed, $\{U\cap\{x_n:n<\omega\}: U\in T\}$ includes all subsets of
the infinite set $\{x_n: n<\omega\}$, which is much more than required (in
particular $\oplus$ holds).

Suppose we have defined $(x_n,U_n)$ for $n<m$ and that there is no
$(x_m,U_m)$ satisfying (a)--(d). This means that if 
$x\in U\in T\cap\uni$, $(\forall n<m)(x_n\notin U)$ and
$x\notin\bigcup_{n<m} U_n$ then
$$ |\{V\in T: (\forall n<m)(x_n\notin V)\mbox{ and }x\notin V\}|<\mu.
\leqno(*)_2$$
Let $U^*=\bigcup\{U\in T\cap\uni: (\forall n<m)(x_n\notin U)\}$.  As
$|\Omega|<\mu={\rm cf}\mu$ we get $$|\{V\in T: (\forall n<m)(x_n\notin V)\
\&\ U^*\setminus(V\cup\bigcup_{n<m}U_n)\neq\emptyset\}|<\mu. \leqno(*)_3$$
Suppose that $U^*\notin\uni$. Then, by $(*)_1$, $\Omega\setminus U^*\in\uni$
(as $U^*$ is open). Since (by clause (c)) $\bigcup_{n<m}U_n\in\uni$ we find
an open set $U$ such that $(\forall n<m)(x_n\notin U)$ and 
$$\mu\leq |\{V\in T: V\cap(\bigcup_{n<m}U_n\cup (\Omega\setminus
U^*))=U\cap(\bigcup_{n<m}U_n\cup(\Omega\setminus U^*))\}|$$
(this is possible by (d)). But if $V\cap(\bigcup_{n<m}U_n\cup (\Omega\setminus
U^*))=U\cap(\bigcup_{n<m}U_n\cup(\Omega\setminus U^*))$,
$V\neq U\cup U^*$ then
$U^*\setminus(V\cup\bigcup_{n<m}U_n)\neq\emptyset$, $(\forall n<m)(x_n\notin
V)$. This contradicts to $(*)_3$. Thus $U^*\in\uni$. Hence (by (d))
we have 
$$\mu\leq |\{V\in T: V\setminus U^*\neq\emptyset\ \&\ (\forall
n<m)(x_n\notin V)\}|.\leqno(*)_4$$ 
Since $|B|<\mu$ we find $V_0\in B$ such that $V_0\setminus
U^*\neq\emptyset$, $(\forall n<m)(x_n\notin V_0)$ and $\mu\leq |\{V\in T:
V_0\subseteq V\}|$. The last condition implies that $\Omega\setminus
V_0\notin\uni$ and hence $V_0\in\uni$ (by $(*)_1$). By the definition of
$U^*$ we conclude $V_0\subseteq U^*$ - a contradiction, thus proving
\ref{22} (when $\lambda=\aleph_0$). \hfill$\Box_{\ref{subcase}}$

\begin{proof}
of \ref{22} when $\lambda>\aleph_0$.
\end{proof}
By Theorem~\ref{third} wlog $|\Omega|=|B|=\lambda$. Let
$\uni=\{A\subseteq\Omega: |\{U\cap A: U\in T\}|\leq \lambda\}$, it is an
ideal.  Let $\uni^+={\cal P}(\Omega)\setminus\uni$.

\begin{observ}
\label{twentysix}
It is enough to prove\\
$\otimes_1$\ \ \ \ for every $Y\in \uni^+$ and $n$ we can find a sequence
$\bar{U}=\langle U_{\zeta}:\zeta<\mu_n\rangle$ of open subsets of $\Omega$
such that one of the following occurs:
\begin{description}
\item[(a)\ ] $\bar{U}$ increasing, $Y\cap U_{\zeta+1}\setminus
U_\zeta\in\uni^+$
\item[(b)\ ] $\bar{U}$ decreasing, $Y\cap U_\zeta\setminus
U_{\zeta+1}\in\uni^+$ 
\item[(c)\ ] $Y\cap U_\zeta\setminus\bigcup_{\vare\neq\zeta}
U_\vare\in\uni^+$
\item[(d)\ ] for some $\langle V_\zeta,y_\zeta:\zeta<\mu_n\rangle$ we have
$Y\cap(\bigcap_{\zeta<\mu_n}U_\zeta\setminus
\bigcup_{\zeta<\mu_n}V_\zeta)\in\uni^+$, 
$V_\zeta$'s and $U_\zeta$'s are open, $V_\zeta\subseteq U_\zeta$,
$y_\zeta\in Y$ are pairwise distinct and
$$U_\zeta\cap\{y_\vare:\vare<\mu_n\} = V_\zeta\cap\{y_\vare:\vare<\mu_n\} =
\{y_\vare:\vare\leq\zeta\}\leqno(*)$$ 
\item[(e)\ ] like (d) but
$$U_\zeta\cap\{y_\vare:\vare<\mu_n\} = V_\zeta\cap\{y_\vare:\vare<\mu_n\} =
\{y_\vare:\zeta \leq \vare < \mu_n\}\leqno(*)'$$ 
\item[(f)\ ] like (d) but
$$U_\zeta\cap\{y_\vare:\vare<\mu_n\} = V_\zeta\cap\{y_\vare:\vare<\mu_n\} =
\{y_\zeta\}\leqno(*)''$$ 
\item[(g)\ ] like (d) but
$$U_\zeta\cap\{y_\vare:\vare<\mu_n\} = V_\zeta\cap\{y_\vare:\vare<\mu_n\} =
\{y_\vare: \vare<\mu_n, \vare\neq\zeta\}\leqno(*)'''$$
\item [(h)\ ] there are $V_\zeta, y_\zeta$ for $\zeta<\mu_n$ such
that $V_\zeta\subseteq U_\zeta$ are open, $y_\zeta\in Y$ are pairwise
distinct, $(U_\zeta\setminus
V_\zeta)\cap\bigcap_{\xi\neq\zeta}V_\xi\in\uni^+$ and 
$$U_\zeta\cap\{y_\vare:\vare<\mu_n\} = V_\zeta\cap\{y_\vare:\vare<\mu_n\} =
\{y_\vare:\vare<\zeta\}\leqno(**)$$ 
\item[(i)\ ] like (h) but
$$U_\zeta\cap\{y_\vare:\vare<\mu_n\} = V_\zeta\cap\{y_\vare:\vare<\mu_n\} =
\{y_\vare:\zeta\leq\vare<\mu_n\}\leqno(**)'$$
\item[(j)\ ] like (h) but 
$$U_\zeta\cap\{y_\vare:\vare<\mu_n\} = V_\zeta\cap\{y_\vare:\vare<\mu_n\} =
\{y_\zeta\}\leqno(**)''$$
\item[(k)\ ] like (h) but
$$U_\zeta\cap\{y_\vare:\vare<\mu_n\} = V_\zeta\cap\{y_\vare:\vare<\mu_n\} =
\{y_\vare:\zeta\neq\vare, \vare<\mu_n\}\leqno(**)'''$$ 
\end{description}
\end{observ}

PROOF: First note that if $n<m<\omega$, $Y_1\subseteq Y_0$, $Y_1,
Y_0\in\uni^+$ and one of the cases (a)--(k) of $\otimes_1$ occurs for $Y_1,
m$ then the same case holds for $Y_0, n$. Consequently, $\otimes_1$ implies
that for each $Y\in\uni^+$ one of (a)--(k) occurs for $Y,n$ for every
$n\in\omega$. Moreover, if $\otimes_1$ then for some
$x\in\{a,b,c,d,e,f,g,h,i,j,k\}$ and $Y_0\in\uni^+$ we have\\ 
$(*)$\ \ \ for every $Y_1\subseteq Y_0$ from $\uni^+$ and $n\in\omega$ case
(x) holds. 

If $x=a$, clause (a) of \ref{22}(1) holds. For this we inductively define
open sets $V_\eta$, $V_\eta^-$ for
$\eta\in\bigcup_{n\in\omega}\prod_{l<n}\mu_l$ such that for
$\eta\in\prod_{l<n}$, $\zeta<\mu_n$: 
\begin{enumerate}
\item $V_\eta^-\subseteq V_\eta$, $(V_\eta\setminus V_\eta^-)\cap
Y_0\in\uni^+$,  $(V_{\eta\hat{\ }\langle \zeta+1\rangle}\setminus
V_{\eta\hat{\ }\langle \zeta\rangle})\cap Y_0\in\uni^+$
\item if $\xi<\mu_{n+1}$ then $V_\eta\subseteq V_{\eta\hat{\
}\langle\zeta\rangle}\subseteq V_{\eta\hat{\
}\langle\zeta,\xi\rangle}\subseteq V^-_{\eta\hat{\ }\langle\zeta+1\rangle}$. 
\end{enumerate}
Let $\langle U_\zeta:\zeta<\mu_0\rangle$ be the increasing sequence of open
sets given by (a) for $Y_0, n=0$. Put
$V_{\langle\zeta\rangle}=U_{2\zeta+1}$,
$V^-_{\langle\zeta\rangle}=U_{2\zeta}$ for $\zeta<\mu_0$. Suppose we have
defined $V_\eta, V^-_\eta$ for $\lg(\eta)\leq m$. Given
$\eta\in\prod_{l<m-1}\mu_l$, $\zeta<\mu_{m-1}$. Apply (a) for
$(V^-_{\eta\hat{\ }\langle\zeta+1\rangle}\setminus V_{\eta\hat{\
}\langle\zeta\rangle}) \cap Y_0$ and $n=m$ to get a sequence $\langle
U_\xi:\xi<\mu_m\rangle$. Put
\[V_{\eta\hat{\ }\langle\zeta,\xi\rangle}=(U_{2\xi+1}\cap V^-_{\eta\hat{\
}\langle\zeta+1\rangle})\cup V_{\eta\hat{\ }\langle\zeta\rangle},\]
\[V^-_{\eta\hat{\ }\langle\zeta,\xi\rangle}=(U_{2\xi}\cap V^-_{\eta\hat{\
}\langle\zeta+1\rangle})\cup V_{\eta\hat{\ }\langle\zeta\rangle}.\]
Next for each $\eta\hat{\
}\langle\zeta\rangle\in\bigcup_{n\in\omega}\prod_{l<n}\mu_l$ choose
$x_{\eta\hat{\ }\langle\zeta\rangle}\in(V_{\eta\hat{\
}\langle\zeta+1\rangle}\setminus V^-_{\eta\hat{\
}\langle\zeta+1\rangle})\cap Y_0$. As the last sets are pairwise disjoint we
get that $x_\eta$'s are pairwise distinct. Moreover, if we put
$U_\eta=\bigcup_{n\in\omega}V_{\eta\rest n}$ (for
$\eta\in\prod_{n\in\omega}\mu_l$)  then we have
\[U_\eta\cap\{x_\nu:\nu\in\bigcup_{n\in\omega}\prod_{l<n}\mu_l\}=\{x_\nu:
(\exists n<\lg(\nu))(\nu\rest n=\eta\rest n\ \&\ \nu(n)<\eta(n))\}.\] 

Similarly one can show that if $x=b$, clause (b) of \ref{22}(1) holds and if
$x=c$ then we can get a discrete set of cardinality $\lambda$ hence all
clauses \ref{22}(1) hold. 

Suppose now that $x=d$. By the induction on $n$ we choose $Y_n$, $\langle
U_{n,\zeta}, V_{n,\zeta}, y_{n,\zeta}:\zeta<\mu_n\rangle$:
\begin{quotation}
\noindent $Y_0=Y$ ($\in\uni^+$)

\noindent $U_{n,\zeta}$, $V_{n,\zeta}$, $y_{n,\zeta}$ (for $\zeta<\mu_n$)
are given by (d) for $Y_n$,

\noindent $Y_{n+1} = Y_n\cap \bigcap_{\zeta<\mu_n} U_{n,\zeta}\setminus 
\bigcup_{\zeta<\mu_n} V_{n,\zeta}\in\uni^+$.
\end{quotation}
For $\eta\in\prod_{l\leq n}\mu_l$ ($n\in\omega$) we let 
\[W_\eta'=V_{n,\eta(n)}\cap\bigcap_{m<n}U_{m,\eta(m)}.\]
As $V_{n,\eta(n)}\cap\{y_{n,\zeta}:\zeta<\mu_n\} =
\{y_{n,\zeta}:\zeta\leq\eta(n)\}$ and $\{y_{n,\zeta}:\zeta<\mu_n\}\subseteq
Y_n\subseteq Y_{m+1}\subseteq U_{m,\eta(m)}$ (for $m<n$) we get
\[W_\eta'\cap\{y_{n,\zeta}:\zeta<\mu_n\}=\{y_{n,\zeta}:\zeta\leq\eta(n)\},\]
\[W_\eta'\cap\{y_{m,\zeta}:\zeta<\mu_m\}\subseteq
U_{m,\eta(m)}\cap\{y_{m,\zeta}:\zeta<\mu_m\} \subseteq \{y_{m,\zeta}:
\zeta\leq\eta(m)\}\ \ \mbox{ (for $m<n$).}\]
Now for $\eta\in\prod_{n<\omega}\mu_n$ we define
$W_\eta=\bigcup_{l<\omega}W_{\eta\rest l}'$. Then for each $n$,
$W_\eta\cap\{y_{n,\zeta}:\zeta<\mu_n\} =\{y_{n,\zeta}: \zeta\leq\eta(n)\}$.
By renaming this implies clause (a) of \ref{22}(1). [For
$\eta\in\prod_{l\leq n}\mu_l$ let $x_\eta=y_{n+1,\gamma(\eta)+1}$, where
$\gamma(\eta)=\mu_n^n\times\eta(0)+\mu_n^{n-1}\times\eta(1)
+\mu_n^{n-2}\times\eta(2) +\ldots+\mu_n^1\times\eta(n-1)+\eta(n)$. Note:
$\mu^l_n$ is the $l$-th ordinal power of $\mu_n$. For
$\eta\in\prod_{l<\omega}\mu_l$ let $\bar{\gamma}(\eta)=0\hat{\ }
\gamma(\eta\rest 1)\hat{\ }\gamma(\eta\rest 2)\hat{\ }\ldots$ and let
$U_\eta=W_{\bar{\gamma}(\eta)}$.]  

For $x=e$ we similarly get clause (b) of \ref{22}(1). For $x=f$ we similarly
get a discrete set of cardinality $\lambda$ so all clauses of 22(1) hold.
The case $x=g$ corresponds to the clause (c) of \ref{22}(1). 

Suppose now that $x=h$. By induction on $n$ we define $Y_\eta$, $U_\eta$,
$V_\eta$ and $x_\eta$ for $\eta\in\prod_{l\leq n}\mu_l$:
\medskip

$Y_{\langle\ \rangle}=Y$,

$U_{\eta\hat{\ }\langle\zeta\rangle}, V_{\eta\hat{\ }\langle\zeta\rangle},
x_{\eta\hat{\ }\langle\zeta\rangle}$ are $U_\zeta, V_\zeta, y_\zeta$ given
by the clause (h) for $Y_\eta$, $\mu_{n+1}$,

$Y_{\eta\hat{\ }\langle\zeta\rangle}=(U_{\eta\hat{\
}\langle\zeta\rangle}\setminus V_{\eta\hat{\
}\langle\zeta\rangle})\cap\bigcap_{\xi\neq\zeta}V_{\eta\hat{\
}\langle\xi\rangle}$.
\medskip

For $\eta\in\prod_{n<\omega}\mu_l$ put
$U_\eta'=\bigcup_{l<\omega}V_{\eta\rest l}$. Then 
\[U_\eta'\cap \{x_\nu:\nu\in\bigcup_{n<\omega}\prod_{l<n}\mu_l\} = \{x_\nu:
\nu=\rho\hat{\ }\langle\zeta\rangle\ \&\ [\neg\rho\triangll\eta\mbox{ or
}\rho\triangll\eta\ \&\ \eta(\mbox{lg\/}(\rho))<\zeta ]\}\] 
witnessing case (e) of \ref{22}(1). 

If $x=i$ then we similarly get case (f) and if $x=j$ we get (d). Lastly
$x=k$ implies the case (c) of \ref{22}(1). 
\hfill$\Box_{\ref{twentysix}}$

\begin{claim}
\label{27}
If $\kappa<\lambda$, $\langle Z_\zeta:\zeta<\kappa\rangle$ is a partition of
$\Omega$, then for some countable $w^*\subseteq\kappa$, for every infinite
$w\subseteq w^*$, $\bigcup_{\zeta\in w}Z_\zeta\notin\uni$. 
\end{claim}

PROOF: \ \ Otherwise there are ${\cal P}\subseteq [\kappa]^{\aleph_0}$ and
$\langle T_w: w\in{\cal P}\rangle$, $T_w\subseteq T$, $|T_w|\leq\lambda$
such that for every $w^*\in [\kappa]^{\aleph_0}$ and $U\in T$, for some
$w\subseteq w^*$, $w\in {\cal P}$ and $V\in T_w$ we have
$U\cap(\bigcup_{\zeta\in w} Z_\zeta)=V\cap(\bigcup_{\zeta\in w}Z_\zeta)$.
Let $\{U_\zeta:\zeta<\lambda\}$ list $\bigcup\{T_w:w\in{\cal P}\}$ (note
that since $\kappa<\lambda$ also
$|[\kappa]^{\leq\aleph_0}|=\kappa^{\aleph_0}<\lambda$). We claim that there is
$U\in T$ such that for every $\xi<\lambda$ there are $\alpha,\beta\in\Omega$
for which:
\begin{description}
\item[(a)\ \ \ ] $\alpha\in U\iff\beta\notin U$
\item[(b)\ \ \ ] $(\forall\zeta<\xi)(\alpha\in U_\zeta\iff\beta\in U_\zeta)$
\item[(c)\ \ \ ] $(\forall \vare<\kappa)(\alpha\in Z_\vare\iff\beta\in
Z_\vare)$ 
\end{description}
Indeed, to find such $U$ consider equivalence relations $E_\xi$ (for
$\xi<\lambda$) determined by (b) and (c), i.e. for $\alpha,\beta\in\Omega$:
\begin{quotation}
\noindent $\alpha\ E_\xi\ \beta$ \ \ if and only if

\noindent $(\forall\zeta<\xi)(\alpha\in U_\zeta\iff\beta\in U_\zeta)$ and

\noindent $(\forall\vare<\kappa)(\alpha\in Z_\vare\iff\beta\in Z_\vare)$.
\end{quotation}
The relation $E_\xi$ has $\leq 2^{|\xi|+\kappa}<\lambda$ equivalence
classes. Consequently for each $\xi<\lambda$
\[|\{V\in T: V\mbox{ is a union of }E_\xi\mbox{-equivalence
classes}\}|<\lambda.\] 
As $|T|>\lambda$ we find a nonempty open set $U$ which for no $\xi<\lambda$
is a union of $E_\xi$-equivalence classes. This $U$ is as needed. 

Now let $(\alpha_n,\beta_n)$ be a pair $(\alpha, \beta)$ satisfying (a)--(c)
for $\xi=\lambda_n$ and let $\{\alpha_n,\beta_n\}\subseteq Z_{\zeta_n}$.
Then $w^*=\{\zeta_n:n<\omega\}$, $U$ contradict the choice of ${\cal P}$ and
$\langle T_w:w\in{\cal P}\rangle$.\hfill$\Box_{\ref{27}}$

\begin{proof}
of $\otimes_1$:\ \ 
\end{proof}
For the notational simplicity we assume that $Y=\Omega$. Let
$B=\bigcup_{n<\omega}B_n$, $|B_n|<\lambda$, $\emptyset\in B_0$.\\ 
As in the proof of claim \ref{lem:an3} wlog for every $x\neq y$ from
$\Omega$ we have
\[|\{U\in T: x\in U\iff y\notin U\}|>\lambda.\]
Let $y_\zeta\in \Omega$ for $\zeta<\mu_{n+6}$ be pairwise distinct. For each
$\zeta <\xi<\mu_{n+6}$ there is $\vare=\vare(\zeta,\xi)\in\{\zeta,\xi\}$ such
that $T^0_{\zeta,\xi}\eqdf\{U\in T: \{y_\zeta,y_\xi\}\cap U=\{y_\vare\}\}$
has cardinality $>\lambda$. For each $U\in T^0_{\zeta,\xi}$ there is
$V[U]\in B$, $y_\vare\in V[U]\subseteq U$. As $|B|\leq\lambda$ for some
$V_{\zeta,\xi}^*\in B$ we have that the set
\[T^1_{\zeta,\xi}=\{U\in T: \{y_\zeta,y_\xi\}\cap U =\{y_\vare\}\mbox{\ \
and\ \ } y_\vare\in V_{\zeta,\xi}^*\subseteq U\}\] 
has cardinality $>\lambda$. For $U\in T$ let $f_U$, $g_U$ be functions such
that:
\begin{enumerate}
\item $f_U:\mu_{n+6}\longrightarrow\omega$, $g_U:\mu_{n+6}\longrightarrow
B$, 
\item $g_U(\vare)=0$ iff $y_\vare\notin U$,
\item if $y_\vare\in U$ then $y_\vare\in g_U(\vare)\subseteq U$,
\item $f_U(\vare)=\min\{n\in\omega:g_U(\vare)\in B_n\}$.
\end{enumerate}
For each $\zeta<\xi<\mu_{n+6}$ we find
$f_{\zeta,\xi}:\mu_{n+6}\longrightarrow\omega$ such that the set 
\[T^2_{\zeta,\xi}=\{U\in T^1_{\zeta,\xi}: f_U=f_{\zeta,\xi}\}\]
has the cardinality $>\lambda$. By Erd\"os-Rado theorem we may assume that
for each $\zeta<\xi<\mu_{n+5}$, $\vare<\mu_{n+5}$ the value of
$f_{\zeta,\xi}(\vare)$ depends on relations between $\zeta,\xi$ and $\vare$
only. Consequently for some $n^*<\omega$, if $\vare<\mu_{n+5}$, $U\in
T^2_{\zeta,\xi}$, $\zeta<\xi<\mu_{n+5}$ then $g_U(\vare)\in B_{n^*}$. As
$|B_{n^*}|<\lambda$ we find (for each $\zeta<\xi<\mu_{n+5}$) a function
$g_{\zeta,\xi}:\mu_{n+6}\longrightarrow B_{n^*}$ such that the set
\[T^3_{\zeta,\xi}=\{U\in T^2_{\zeta,\xi}: g_U=g_{\zeta,\xi}\}\]
is of the size $>\lambda$. Let 
\[U_{\zeta,\xi}=\bigcup T^3_{\zeta,\xi},\ \ \ \ \
V_{\zeta,\xi}=\bigcup_{\vare<\mu_{n+5}}g_{\zeta,\xi}(\vare).\] 
Clearly 
$$V_{\zeta,\xi}\subseteq U_{\zeta,\xi}, U_{\zeta,\xi}\setminus
V_{\zeta,\xi}\notin\uni, U_{\zeta,\xi}\cap\{y_\zeta,y_\xi\}=
V_{\zeta,\xi}\cap\{y_\zeta,y_\xi\}=\{y_\vare(\zeta,\xi)\}\mbox{
and}\leqno(*)$$  
$$U_{\zeta,\xi}\cap\{y_\delta:\delta<\mu_{n+5}\}=
V_{\zeta,\xi}\cap\{y_\delta:\delta<\mu_{n+5}\}.\leqno(**)$$
Let $T_1=\{V_{\zeta,\xi},U_{\zeta,\xi}: \zeta<\xi<\mu_{n+5}\}$, so
$|T_1|<\lambda$. Define a two place relation $E_{T_1}$ on $\Omega$:
\[x\/E_{T_1}\/y\ \mbox{ iff }\ (\forall U\in T_1)(x\in U\iff y\in U).\] 
Clearly $E_{T_1}$ is an equivalence relation with $\leq
2^{|T_1|}<\lambda$ equivalence classes.
Hence by claim~\ref{27} for each $\zeta<\xi<\mu_{n+5}$, for some
$\omega$-sequence of 
$E_{T_1}$-equivalence classes $\langle A_{\zeta,\xi,n}:n<\omega\rangle$ we
have:
\[A_{\zeta,\xi,n}\subseteq U_{\zeta,\xi}\setminus V_{\zeta,\xi}\mbox{\ \
and\ \ for each infinite } w\subseteq\omega, \bigcup_{n\in w}
A_{\zeta,\xi,n}\notin\uni.\] 
By Erd\"os-Rado theorem, wlog for $\zeta_1<\zeta_2<\mu_{n+4}$,
$\xi_1,\xi_2<\mu_{n+4}$ the truth values of
``$\vare(\zeta_1,\zeta_2)=\zeta_1$'', ``$y_{\xi_1}\in
V_{\zeta_1,\zeta_2}$'', ``$y_{\xi_1}\in U_{\zeta_1,\zeta_2}$'',
``$A_{\zeta_1,\zeta_2,n}\subseteq U_{\xi_1,\xi_2}$'',
``$A_{\zeta_1,\zeta_2,n}\subseteq V_{\xi_1,\xi_2}$'', 
``$A_{\zeta_1,\zeta_2,n}=A_{\zeta_1,\zeta_2,m}$'',
``$A_{\zeta_1,\zeta_2,n}=A_{\xi_1,\xi_2,m}$'' depend just on the order and
equalities among $\zeta_1,\zeta_2,\xi_1,\xi_2$ (and of course $n$, $m$).

As each infinite union $\bigcup_{n\in\omega}A_{\zeta,\xi,n}$ is large, wlog
those truth values also does not depend on $n$ (for the last one we mean
``$A_{\zeta_1,\zeta_2,n}=A_{\xi_1,\xi_2,n}$''). Note: if
$A_{1,2,n}=A_{3,4,m}$ then $A_{1,2,n}=A_{3,4,n}=A_{1,2,m}$.

Now, $A_{\zeta,\xi,n}$ is either included in $U_{\xi_1,\xi_2}$ or is
disjoint from it (uniformly for $n$); similarly for $V_{\xi_1,\xi_2}$.

\underline{Case A}: \ \ $A_{3,4,n}\cap U_{1,2}=\emptyset$

\noindent Let $U_\zeta'=\bigcup_{\xi\leq\zeta}U_{2\xi,2\xi+1}$. Then
$\langle U_\zeta':\zeta<\mu_n\rangle$ is an increasing sequence of open sets
and $\bigcup_{n\in\omega}A_{2\zeta+2,2\zeta+3,n}\subseteq
U_{\zeta+1}'\setminus U_\zeta'$, which witnesses that the last set is in
$\uni^+$. Thus we get clause (a).

\underline{Case B}: \ \ $A_{1,2,n}\cap U_{3,4}=\emptyset$

\noindent Let $U_\zeta'=\bigcup_{\zeta\leq\xi<\mu_n}U_{2\xi,2\xi+1}$. Then
$\langle U_\zeta':\zeta<\mu_n\rangle$ is a decreasing sequence of open sets
and $\bigcup_{n\in\omega}A_{2\zeta,2\zeta+1,n}\subseteq U_\zeta'\setminus
U_{\zeta+1}'$. Consequently we get clause (b).

Thus we have to consider the case
$$A_{1,2,n}\subseteq U_{3,4}\ \mbox{ and }\ A_{3,4,n}\subseteq
U_{1,2}$$ 
only. So we assume this.

\underline{Case C}: \ \ $A_{1,2,n}\cap V_{3,4}=\emptyset$, $A_{3,4,n}\cap
V_{1,2}=\emptyset$ 

Let $U_\zeta'=U_{2\zeta,2\zeta+1}$, $V_\zeta'=V_{2\zeta,2\zeta+1}$.
\begin{description}
\item[subcase C1:]\ \ \ $y_1\in U_{3,4}$, $y_5\in U_{3,4}$\\
Then let $y_\zeta'$ is the unique member of
$\{y_{2\zeta},y_{2\zeta+1}\}\setminus\{y_{\vare(2\zeta,2\zeta+1)}\}$.\\ 
By $(**)$ we easily get that $\langle U_\zeta',V_\zeta',y_\zeta':
\zeta<\mu_n\rangle$ witnesses the clause (g).
\item[subcase C2:]\ \ \ either $y_1\notin U_{3,4}$ or $y_5\notin U_{3,4}$\\
Then we put $y_\zeta'=y_{\vare(2\zeta,2\zeta+1)}$ and we get one of the
cases (d), (e) or (f).
\end{description}

\underline{Case D}:\ \ $A_{1,2,n}\subseteq V_{3,4}$, $A_{3,4,n}\cap
V_{1,2}=\emptyset$ 

\noindent We let $U_\zeta'=\bigcup\{V_{2\xi,2\xi+1}:\xi\leq\zeta\}$. Thus
$U_\zeta'$ increases with $\zeta$ and $U_{\zeta+1}'\setminus U_\zeta'$
includes $\bigcup_{n\in\omega}A_{2\zeta,2\zeta+1,n}$. Thus clause (a) holds.

\underline{Case E}:\ \ $A_{1,2,n}\cap V_{3,4}=\emptyset$, $A_{3,4,n}\subseteq
V_{1,2}$ 

\noindent Let $U_\zeta'=\bigcup\{V_{2\xi,2\xi+1}:\xi\geq\zeta\}$. Then
$U_\zeta'$ decrease with $\zeta$ and the clause (b) holds. 

\underline{Case F}:\ \ $A_{1,2,n}\subseteq V_{3,4}$, $A_{3,4,n}\subseteq
V_{1,2}$ 

\noindent Let $U_\zeta'=U_{2\zeta,2\zeta+1}$,
$V_\zeta'=V_{2\zeta,2\zeta+1}$. If $y_1,y_5\in U_{3,4}$ then we put
$y_\zeta'\in\{y_{2\zeta},y_{2\zeta+1}\}\setminus\{y_{\vare(2\zeta,2\zeta+1)}\}$
and we get case (k). Otherwise we put $y_\zeta'=y_{\vare(2\zeta,2\zeta+1)}$
and we obtain one of the cases (h), (i) or (j). \hfill$\Box_{\ref{22}}$ 

\begin{concluding}
\begin{enumerate}
\item Assume that a topology $T$ on $\Omega$ with a base $B$ and $\lambda$,
$\langle\mu_n:n\in\omega\rangle$ are as before ($\mu_n$ regular for
simplicity). If 
\begin{description}
\item[(*)] $x_\nu\in\Omega$ for $\nu\in\bigcup_{n\in\omega}\prod_{l<n}\mu_l$
and $U_\eta\in T$ for $\eta\in\prod_{n\in\omega}\mu_n$ and
\item[(**)] if $n<\omega$, $\nu\in\prod_{l<n}\mu_l$ and
$\eta\in\prod_{l<\omega}\mu_l$ then for some $k$, 
\[(\forall\eta')(\eta'\in\prod_{l<\omega}\mu_l\ \&\ \eta'\restriction
k=\eta\restriction k\ \Rightarrow\
U_{\eta'}\cap\{x_\nu\}=U_\eta\cap\{x_\nu\}).\] 
\end{description}
Then we can find $S\subseteq\bigcup_{n<\omega}\prod_{l<n}\mu_l$ and $\langle
U_{\eta,\nu}: \eta,\nu\in\prod_{l<n}\mu_l\cap S\mbox{\ \ \ for some }
n\rangle$ and $\langle U^*_\eta:\eta\in\lim S\rangle$ (where $\lim
S=\{\eta\in\prod_{l<\omega}\mu_l: (\forall l<\omega)(\eta\restriction l\in
S)\}$) such that 
\begin{description}
\item [(a)] $\langle\rangle\in S$, $S$ is closed under initial segments and 
\[\eta\in S\ \&\ n={\rm lg}\eta\ \Rightarrow\ (\exists\alpha)(\eta\hat{\
}\langle \alpha\rangle\in S)\]
and for some infinite $w\subseteq\omega$, for every $n<\omega$ and
$\eta\in\lim S$ we have:
\[n\in w\iff (\exists^{\geq 2}\alpha<\mu_n)(\eta\hat{\ }\langle\alpha\rangle
\in S) \iff (\exists^{\mu_n}\alpha<\mu_n)(\eta\hat{\
}\langle\alpha\rangle\in S).\]
\item[(b)] if $\rho,\nu\in\prod_{l<n}\mu_l\cap S$ and $\nu\triangll\eta\in
S\cap\prod_{l<\omega}\mu_l$ then $U^*_\eta\cap\{x_\rho\}=
U^*_{\nu,\eta}\cap\{x_\rho\}$,
\item[(c)] for $\eta\in\lim S$, $U^*_\eta\cap\{x_\rho:\rho\in S\}=
U_\eta\cap\{x_\rho: \rho\in S\}$.
\end{description}
\item So in Theorem~\ref{22}, the case (c) can be further described. 
\item We can consider basic forms for any analytic families of subsets of
$\lambda$ (then we have more cases; as in 23 and $\otimes_1$ of 26).
\end{enumerate}
\end{concluding}

\shlhetal

\end{document}